\documentclass[12pt]{amsart}

\usepackage{latexsym, amssymb,lscape,graphicx,epsfig}

\theoremstyle{plain}
\newtheorem{theorem}{Theorem}[section]
\newtheorem{proposition}{Proposition}[section]

\theoremstyle{definition}
\newtheorem{definition}{Definition}[section]
\newtheorem{lem}{Lemma}[section]

\newtheorem{cor}{Corollary}
\newtheorem{example}{Example}
\newtheorem{exercise}{Exercise}

\newenvironment{theo}[1]%
{\begin{theorem}\label{T:#1}}%
{\end{theorem}}

{\begin{proposition}\label{T:#1}}%
{\end{proposition}}

{\vskip6pt}

\newenvironment{remark}%
{\vskip6pt%
\noindent%
{\it Remark.}}%
{\vskip6pt}

{\vskip6pt%
\noindent%
{\it Remarks.}}%
{\vskip6pt}

{\begin{exercise}\label{T:#1}}%
{\end{exercise}}

\newenvironment{ex}[1]%
{\begin{example}\label{T:#1}}%
{\end{example}}

\newcommand{\RR}{\text{${\mathbb R}$}}

\newcommand{\Ll}{\mathcal{L}}

\newcommand{\C}{\text{$\mathbb C$}}

\newcommand{\ZZ}{\mathbb Z}

\renewcommand{\frak}[1]{\text{$\mathfrak{#1}$}}

\newcommand{\ga}{\text{$\alpha$}}
\newcommand{\gb}{\text{$\beta$}}

\newcommand{\gO}{\text{$\Omega$}}
\newcommand{\go}{\text{$\omega$}}
\newcommand{\om}{\text{$\omega$}}

\newcommand{\Gg}{\mathfrak{g}}

\newcommand{\del}{\text{$\partial$}}
\newcommand{\delbar}{\overline{\partial}}

\newcommand{\Ann}{\mathrm{Ann}}

\newcommand{\nil}{\mathrm{nil}}

\newcommand{\mc}[1]{\text{$\mathcal{#1}$}}

\newcommand{\gcss}{generalized complex structures}

\newcommand{\into}{\rightarrow}
\newcommand{\IP}[1]{\langle #1\rangle}
\newcommand{\JJ}{\mathcal{J}}
\newcommand{\End}{\text{End}}
\newcommand{\idx}{\mathrm{nil}}

\numberwithin{equation}{section}
\author{Gil R. Cavalcanti \and Marco Gualtieri}
\email{gilrc@maths.ox.ac.uk}
\email{mgualtie@fields.utoronto.ca}
\title{Generalized complex structures on nilmanifolds}
\begin{document}
\begin{abstract}
We show that all 6-dimensional nilmanifolds admit generalized
complex structures.  This includes the five classes of nilmanifold
which admit no known complex or symplectic structure. Furthermore,
we classify all 6-dimensional nilmanifolds according to which of
the four types of left-invariant generalized complex structure
they admit.  We also show that the two components of the
left-invariant complex moduli space for the Iwasawa manifold are
no longer disjoint when they are viewed in the generalized complex
moduli space. Finally, we provide an 8-dimensional nilmanifold
which admits no left-invariant generalized complex structure.
\end{abstract}
\maketitle

\markboth{\sc G. R. Cavalcanti and M. Gualtieri}{\sc Generalized
complex structures on nilmanifolds}

\section*{Introduction}

Ever since Thurston~\cite{Th} presented a nilmanifold as the
first instance of a symplectic but non-K\"ahler manifold in 1976,
the study of invariant geometries on nilmanifolds has been an
interesting source of examples in differential geometry.

A nilmanifold is a homogeneous space $M=\Gamma\backslash G$, where
$G$ is a simply connected nilpotent real Lie group and $\Gamma$ is
a lattice of maximal rank in $G$.
Such groups $G$ of dimension $\leq 7$ have been classified, and 6
is the highest dimension where there are finitely many. According
to \cite{Mag,Mo} there are 34 isomorphism classes of connected,
simply-connected 6-dimensional nilpotent Lie groups. This means
that, with respect to invariant geometry, there are essentially 34
separate cases to investigate.

The question of which 6-dimensional nilmanifolds admit
 symplectic structure was settled by Goze and
Khakimdjanov \cite{GK}: exactly 26 of the 34 classes admit
 symplectic forms.  Subsequently, the question of
left-invariant complex geometry was solved by Salamon~\cite{Sa}:
he proved that exactly 18 of the 34 classes admit invariant
complex structure.  While the torus is the only nilmanifold
admitting K\"ahler structure, 15 of the 34 nilmanifolds admit both
complex and symplectic structures.  This leaves us with 5 classes
of 6-dimensional nilmanifolds admitting neither complex nor
symplectic left-invariant geometry.  See Figure 1 for
illustration.

It was this result of Salamon which inspired us to ask whether the
5 outlying classes might admit \emph{generalized complex
structure}, a geometry recently introduced by Hitchin~\cite{Hi}
and developed by the second author~\cite{Marco}.  Generalized
complex geometry contains complex and symplectic geometry as
extremal special cases and shares important properties with them,
such as an elliptic deformation theory as well as a local normal
form (in regular neighbourhoods).  The main result of this paper
is to answer this question in the affirmative: all 6-dimensional
nilmanifolds admit generalized complex structures.

\setlength{\unitlength}{1mm}
\begin{picture}(60,51)(-65,-26)
\linethickness{1pt}
\thinlines
\qbezier(10,15)(25,15)(25,0)
\qbezier(25,0)(25,-15)(10,-15)
\qbezier(10,-15)(-5,-15)(-5,0)
\qbezier(-5,0)(-5,15)(10,15)
\qbezier(-10,15)(-25,15)(-25,0)
\qbezier(-25,0)(-25,-15)(-10,-15)
\qbezier(-10,-15)(5,-15)(5,0)
\qbezier(5,0)(5,15)(-10,15)
\put(13,10){\makebox(0,0){\footnotesize Symplectic}}
\put(15,0){\makebox(0,0){$11$}}
\put(-15,0){\makebox(0,0){$3$}}
\put(0,0){\makebox(0,0){$15$}}
\put(0,-18){\makebox(0,0){$5$}}
\put(-11,10){\makebox(0,0){\footnotesize Complex}}
\put(-30,-25){\framebox(60,45)[bc]{{\footnotesize Generalized complex}}}
\end{picture}
\begin{center}
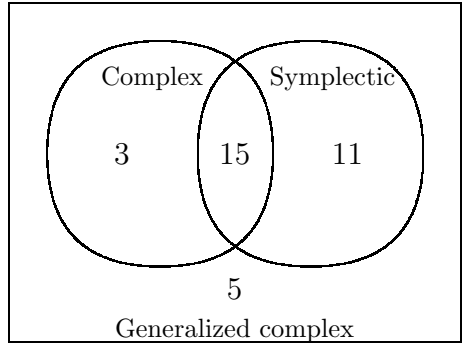
{\small {\it Figure 1: Left-invariant structures on the 34
six-dimensional nilpotent Lie groups.}}
\end{center}
\vskip12pt

We begin in Section~\ref{gcs} with a review of generalized complex
geometry.  A brief introduction to nilmanifolds follows in
Section~\ref{nil}. Some results about generalized complex
structures on nilmanifolds in arbitrary dimension appear in
Section~\ref{gen}. Section~\ref{six} contains our main result: the
classification of left-invariant generalized complex structures on
6-dimensional nilmanifolds. In Section~\ref{iwa} we show that
while the moduli space of left-invariant complex structures on the
Iwasawa nilmanifold is disconnected (as shown in~\cite{KS}), its
components can be joined using generalized complex structures. In
the final section, we provide an 8-dimensional nilmanifold which
does not admit a left-invariant generalized complex structure,
thus precluding the possibility that all nilmanifolds admit
left-invariant generalized complex geometry.

The authors wish to thank Nigel Hitchin and Simon Salamon for
helpful discussions.  The first author was supported by CAPES
grant 1326/99-6, and the second by an NSERC fellowship.
%

\section{Generalized complex structures}\label{gcs}

We briefly review the theory of generalized complex structures;
see~\cite{Marco} for details. A generalized complex structure on a
smooth manifold $M$ is defined to be a complex structure $\JJ$,
not on the tangent bundle $T$, but on the sum $T\oplus T^*$ of the
tangent and cotangent bundles. This complex structure is required
to be orthogonal with respect to the natural inner product on
sections $X+\xi,Y+\eta\in C^\infty(T\oplus T^*)$ defined by
\begin{equation*}
\IP{X+\xi,Y+\eta}=\tfrac{1}{2}(\xi(Y)+\eta(X)).
\end{equation*}
This is only possible if the manifold has even dimension, so we
suppose $\dim_{\RR} M = 2n$. In addition, the $+i$-eigenbundle
\begin{equation*}
L < (T\oplus T^*)\otimes\C
\end{equation*}
of $\JJ$ is required to be involutive with respect to the Courant
bracket, a skew bracket operation on smooth sections of $T\oplus
T^*$ defined by
\begin{equation*}
[X+\xi,Y+\eta]=[X,Y]+\Ll_X\eta-\Ll_Y\xi -
\tfrac{1}{2}d(i_X\eta-i_Y\xi),
\end{equation*}
where $\Ll_X$ and $i_X$ denote the Lie derivative and interior
product operations on forms.

Since $\JJ$ is orthogonal with respect to $\IP{\cdot,\cdot}$, the
$+i$-eigenbundle $L$ is a maximal isotropic sub-bundle of
$(T\oplus T^*)\otimes\C$, and as such can be expressed as the
Clifford annihilator of a unique line sub-bundle $U_L$ of the
complex spinors for the metric bundle $T\oplus T^*$.  Since its
annihilator is maximal isotropic, $U_L$ is by definition a
\emph{pure} spinor line, and we call it the \emph{canonical} line
bundle of $\JJ$.

The bundle $\wedge^\bullet T^*$ of differential forms can in fact
be viewed as a spinor bundle for $T\oplus T^*$, where the Clifford
action of an element $X+\xi\in T\oplus T^*$ on a differential form
$\rho$ is given by
\begin{equation*}
(X+\xi)\cdot \rho = i_X\rho + \xi\wedge\rho.
\end{equation*}
Note that $(X+\xi)^2\cdot\rho = \IP{X+\xi,X+\xi}\rho$, as
required.  Therefore, the canonical bundle $U_L$ may be viewed as
a smooth line sub-bundle of the complex differential forms
according to the relation
\begin{equation}\label{relate}
L=\left\{X+\xi\in (T\oplus T^*)\otimes\C\ :\ (X+\xi)\cdot U_L=0
\right\}.
\end{equation}

At every point, the line $U_L$ is generated by a complex
differential form of special algebraic type: purity is equivalent
to the fact that it has the form
\begin{equation}\label{algebraic}
\rho=e^{B+i\omega}\Omega,
\end{equation}
where $B,\omega$ are real $2$-forms and
$\Omega=\theta_1\wedge\cdots\wedge\theta_k$ is a complex
decomposable $k$-form. The number $k$ is called the \emph{type} of
the generalized complex structure, and it is not required to be
constant along the manifold.  Points where the type is locally
constant are called \emph{regular}.  Since $L$ is the
$+i$-eigenbundle of a complex structure, we see that
$L\cap\overline{L}=\{0\}$. This is equivalent to an additional
constraint on $\rho$:
\begin{equation}\label{nondegeneracy}
\omega^{2n-2k}\wedge\Omega\wedge\overline\Omega\neq 0.
\end{equation}
Hence we see that on a $2n$-manifold the type may take values from
$k=0$ to $k=n$.  Finally, as is proven in~\cite{Marco}, the
involutivity of $L$ under the Courant bracket is equivalent to the
condition, on any local trivialization $\rho$ of $U_L$, that there
exist a section $X+\xi\in C^\infty(T\oplus T^*)$ such that
\begin{equation}\label{integrability}
d\rho = (X+\xi)\cdot\rho.
\end{equation}
Near any regular point, this condition implies that the
distribution determined by $\Omega\wedge\overline{\Omega}$
integrates to a foliation, and with (\ref{nondegeneracy}), also
implies that $\omega$ is a leafwise symplectic form.

In the special case that $U_L$ is a trivial bundle admitting a
nowhere-vanishing \emph{closed} section $\rho$, the structure is
said to be a \emph{generalized Calabi-Yau} structure, as
in~\cite{Hi}.

\subsection{Examples}
So far, we have explained how a generalized complex structure is
equivalent to the specification of a pure line sub-bundle of the
complex differential forms, satisfying the non-degeneracy
condition~(\ref{nondegeneracy}) and the integrability
condition~(\ref{integrability}).  Now let us provide some examples
of such structures.

\begin{example}[Complex geometry (type $n$)]
Let $J\in\End(T)$ be a usual complex structure on a $2n$-manifold.
The generalized complex structure corresponding to $J$ is
\begin{equation*}
\JJ_J=\left(\begin{matrix}-J&0\\0&J^*\end{matrix}\right),
\end{equation*}
where the matrix is written in the natural splitting $T\oplus
T^*$. Clearly $\JJ_J^2=-1$, and orthogonality is easily verified.
The $+i$-eigenbundle of $\JJ_J$ is the maximal isotropic
\begin{equation*}
L=T_{0,1}\oplus T^*_{1,0},
\end{equation*}
where $T_{1,0}=\overline{T_{0,1}}$ is the $+i$-eigenbundle of $J$
in the usual way. The bundle $L$ is the Clifford annihilator of the
line bundle
\begin{equation*}
U_L=\wedge^n(T^*_{1,0}),
\end{equation*}
the canonical bundle associated to $J$.  We see that $\JJ_J$ is of
type $n$ at all points in the manifold. The Courant involutivity
of $L$ is equivalent to the Lie involutivity of $T_{0,1}$, which
is the usual integrability condition for complex structures.  To
be generalized Calabi-Yau, there must be a closed trivialization
$\Omega\in C^\infty(U_L)$, which is the usual Calabi-Yau
condition.
\end{example}

\begin{example}[Symplectic geometry (type $k=0$)]
Let $\omega\in \Omega^2(M)$ be a usual symplectic structure,
viewed as a skew-symmetric isomorphism $\omega:T\rightarrow T^*$
via the interior product $X\mapsto i_X\omega$. The generalized
complex structure corresponding to $\omega$ is
\begin{equation*}
\JJ_\omega=\left(\begin{matrix}0&-\omega^{-1}\\\omega&0\end{matrix}\right),
\end{equation*}
where the matrix is written in the natural splitting $T\oplus
T^*$. Clearly $\JJ_\omega^2=-1$, and orthogonality is easily
verified. The $+i$-eigenbundle of $\JJ_\omega$ is the maximal
isotropic
\begin{equation*}
L=\left\{X-i\omega(X)\ :\ X\in T\otimes\C\right\},
\end{equation*}
which is the Clifford annihilator of the line bundle $U_L$ with
trivialization given by
\begin{equation*}
\rho = e^{i\omega}.
\end{equation*}
We see that $\JJ_\omega$ is of type $0$ at all points in the
manifold.  The Courant involutivity of $L$ is equivalent to the
constraint $d\rho=0$, itself equivalent to the usual integrability
condition $d\omega=0$ for symplectic structures.  Note that
symplectic structures are always generalized Calabi-Yau.
\end{example}

The preceding examples demonstrate how complex and symplectic
geometry appear as extremal cases of generalized complex geometry.
We will now explain how one may deform these examples into new
ones.

\subsection{$B$-fields and $\beta$-fields}
Unlike the Lie bracket, whose only symmetries are diffeomorphims,
the Courant bracket is preserved by an additional group of
symmetries of $T\oplus T^*$, consisting of closed $2$-forms $B$
acting via the orthogonal shear transformation
\begin{equation*}
X+\xi\mapsto X+\xi-i_XB.
\end{equation*}
Such automorphisms are called \emph{B-field} transformations, and
their associated spinorial action on differential forms is via the
exponential:
\begin{equation*}
\rho\mapsto e^B\wedge\rho.
\end{equation*}
In this way, we see that $B$-field transformations do not have any
effect on the type of a generalized complex structure.
Nevertheless, using $B$-fields and diffeomorphisms, one may choose
canonical coordinates for a generalized complex structure, near
any regular point:
\begin{theorem}[\cite{Marco}, Theorem 4.35] Any regular point of type $k$ in
a generalized complex $2n$-manifold has a neighbourhood which is
equivalent, via a diffeomorphism and a $B$-field transformation, to
the product of an open set in $\C^k$ with an open set in the
standard symplectic space $\RR^{2n-2k}$.
\end{theorem}

Although automorphisms of the Courant bracket do not affect the
type of a generalized complex structure, there may be
non-automorphisms which nevertheless transform a given generalized
complex structure into another one, of modified type. For example,
consider the action of a bivector $\beta\in C^\infty(\wedge^2 T)$
on $T\oplus T^*$ via the orthogonal shear transformation
\begin{equation*}
X+\xi\mapsto X+\xi+i_{\xi}\beta.
\end{equation*}
The spinorial action of such $\beta$-\emph{field} transformations
on differential forms is also via the exponential
\begin{equation*}
\rho\mapsto e^{\beta}\rho = (1+i_\beta +
\tfrac{1}{2}i_\beta^2+\cdots)\rho.
\end{equation*}
The following proposition describes the conditions on $\beta$
which ensure that it takes a complex structure into a generalized
complex structure of different type.
\begin{proposition}\label{T:beta field}{\em(\cite{Marco}, Section 5.3)}
Let $J$ be a complex structure on a compact $2n$-manifold, viewed
as a generalized complex structure of type $n$. Let $\beta$ be a
smooth bivector of type $(2,0)$ with respect to $J$. Then by the
above action, $\beta$ deforms $J$ into another generalized complex
structure if and only if it is sufficiently small and
\begin{equation*}
\delbar\beta + \tfrac{1}{2}[\beta,\beta]=0,
\end{equation*}
which holds if and only if each summand vanishes separately, i.e.
$\beta$ is a holomorphic Poisson structure.  The resulting
generalized complex structure has type $n-k$ at points where the
bivector $\beta$ has rank $k$.
\end{proposition}

For example, any holomorphic bivector $\beta$ on $\C P^2$ is
Poisson, and therefore a sufficiently small constant multiple of
it will deform the standard complex structure into a generalized
complex structure. Since $\beta$ vanishes on a cubic and is of
rank 2 elsewhere, the resulting generalized complex structure has
type 2 along the cubic and is of type 0 elsewhere.  Our purpose in
introducing $\beta$-transforms is that we will use them to produce
some examples of generalized complex structures on nilmanifolds.

\section{Nilmanifolds}\label{nil}

A nilmanifold is the quotient $M=\Gamma\backslash G$ of a
connected, simply-connected nilpotent real Lie group $G$ by the
left action of a maximal lattice $\Gamma$, i.e. a discrete
cocompact subgroup. By results of Malcev \cite{Ma}, a nilpotent
Lie group admits such a lattice if and only if its Lie algebra has
rational structure constants in some basis.  Moreover, any two
nilmanifolds of $G$ can be expressed as finite covers of a third
one.

A connected, simply-connected nilpotent Lie group is diffeomorphic
to its Lie algebra via the exponential map and so is contractible.
For this reason, the homotopy groups $\pi_k$ of nilmanifolds
vanish for $k>1$, i.e. nilmanifolds are Eilenberg-MacLane spaces
$K(\Gamma,1)$.  In fact, their diffeomorphism type is determined
by their fundamental group.  Malcev showed that this fundamental
group is a finitely generated nilpotent group with no element of
finite order.  Such groups can be expressed as central $\ZZ$
extensions of groups of the same type, which implies that any
nilmanifold can be expressed as a circle bundle over a nilmanifold
of lower dimension.  Because of this, one may easily use Gysin
sequences to compute the cohomology ring of any nilmanifold.
Nomizu used this fact to show that the rational cohomology of a
nilmanifold is captured by the subcomplex of the de Rham complex
$\Omega^\bullet(M)$ consisting of forms descending from
left-invariant forms on $G$:

\begin{theo}{Nomizu}
{\em (Nomizu \cite{No})} The de Rham complex $\Omega^\bullet(M)$
of a nilmanifold $M=\Gamma\backslash G$ is quasi-isomorphic to the
complex $\wedge^{\bullet}\Gg^*$ of left-invariant forms on $G$,
and hence the de Rham cohomology of $M$ is isomorphic to the Lie
algebra cohomology of~$\Gg$.
\end{theo}

In this paper we will search for generalized complex structures on
$\Gamma\backslash G$ which descend from left-invariant ones on
$G$, which we will call \emph{left-invariant} generalized complex
structures. This will require detailed knowledge of the structure
of the Lie algebra $\Gg$, and so we outline its main properties in
the remainder of this section.

Nilpotency implies that the central descending series of ideals
defined by $\Gg^0=\Gg,\ \Gg^k = [\Gg^{k-1},\Gg]$ reaches $\Gg^s=0$
in a finite number $s$ of steps, called the {\it nilpotency index}, $\nil(\Gg)$  (also called the nilpotency index of any nilmanifold associated to $\Gg$). Dualizing, we obtain a filtration of $\Gg^*$ by the
annihilators $V_i$ of $\Gg^i$, which can also be expressed as
\begin{equation*}
V_i = \left\{v\in\Gg^*\ :\ dv\in V_{i-1}\right\},
\end{equation*}
where $V_0=\{0\}$.  Choosing a basis for $V_1$ and extending
successively to a basis for each $V_k$, we obtain a {\it Malcev
basis} $\{e_1,\ldots,e_{n}\}$ for $\Gg^*$. This basis satisfies
the property
\begin{equation}\label{E:basis}
de_i \in {\wedge}^2\left<e_1, \ldots, e_{i-1}\right>\ \ \forall i.
\end{equation}
The filtration of $\Gg^*$ induces a filtration of its exterior
algebra, and leads to the following useful definition:
\begin{definition}
With $V_i$ as above, the {\it nilpotent degree} of a $p$-form \ga,
which we denote by  $\idx(\ga)$, is the smallest $i$ such that
$\ga \in {\wedge}^p V_i$.
\end{definition}

\begin{remark}
If \ga\ is a 1-form of nilpotent degree $i$ then $\idx(d\ga) =
i-1$.
\end{remark}

In this paper, we specify the structure of a particular nilpotent
Lie algebra by listing the exterior derivatives of the elements of
a Malcev basis as an $n$-tuple of 2-forms $(de_k = \sum c^{ij}_k
e_i\wedge e_j)_{k=1}^m$. In low dimensions we use the shortened
notation $ij$ for the 2-form $e_i\wedge e_j$, as in the following
6-dimensional example: the 6-tuple $(0,0,0,12,13,14+35)$ describes
a nilpotent Lie algebra with dual $\Gg^*$ generated by 1-forms
$e_1,\ldots, e_6$ and such that $de_1=de_2=de_3=0$, while
$de_4=e_1\wedge e_2$, $de_5=e_1\wedge e_3$, and $de_6=e_1\wedge
e_4+e_3\wedge e_5$.  We see clearly that
$V_1=\left<e_1,e_2,e_3\right>$,
$V_2=\left<e_1,e_2,e_3,e_4,e_5\right>$, and $V_3=\Gg^*$, showing
that the nilpotency index of $\Gg$ is 3.

\section{Generalized complex structures on
nilmanifolds}\label{gen}

In this section, we present two results concerning generalized
complex structures on nilmanifolds of arbitrary dimension.  In
Theorem \ref{pure spinor}, we prove that any left-invariant
generalized complex structure on a nilmanifold must be generalized
Calabi-Yau, i.e. the canonical bundle $U_L$ has a closed
trivialization.  In Theorem \ref{T:maximal} we prove an upper bound
for the type of a left-invariant generalized complex structure,
depending only on crude information concerning the nilpotent
structure.

We begin by observing that a left-invariant generalized complex
structure must have constant type $k$ throughout the nilmanifold
$M^{2n}$, and its canonical bundle $U_L$ must be trivial.  Hence,
by (\ref{algebraic}),(\ref{nondegeneracy}), and
(\ref{integrability}) we may choose a global trivialization of the
form
\begin{equation*}
\rho = e^{B+i\om}\Omega,
\end{equation*}
where $B,\omega$ are real left-invariant 2-forms and $\Omega$ is a
globally decomposable complex $k$-form, i.e.
\begin{equation*}
\Omega=\theta_1\wedge\cdots\wedge\theta_k,
\end{equation*}
with each $\theta_i$ left-invariant. These data satisfy the
nondegeneracy condition
$\omega^{2n-2k}\wedge\Omega\wedge\overline{\Omega}\neq 0$ as well
as the integrability condition $d\rho = (X + \xi)\cdot\rho$ for
some section $X + \xi\in C^\infty(T\oplus T^*)$.  Since $\rho$ and
$d\rho$ are left-invariant, we can choose $X+\xi$ to be
left-invariant as well.

It is useful to order $\{\theta_1,\ldots,\theta_k\}$ according to
nilpotent degree, and also to choose them in such a way that
$\{\theta_j\ : \ \idx(\theta_j)>i\}$ is linearly independent
modulo $V_i$; this is possible according to the following lemma.
\begin{lem}\label{T:trivial}
It is possible to choose a left-invariant decomposition
$\Omega=\theta_1\wedge\cdots\wedge\theta_k$ such that
$\idx(\theta_i)\leq\idx(\theta_j)$ for $i<j$, and such that
$\{\theta_j\ :\ \idx(\theta_j) > i\}$ is linearly independent
modulo $V_i$ for each $i$.
\end{lem}
\begin{proof}
Choose any left-invariant decomposition
$\Omega=\theta_1\wedge\cdots\wedge\theta_k$ ordered according to
nilpotent degree, i.e. $\idx(\theta_i)\leq\idx(\theta_j)$ for
$i<j$.  Certainly $\{\theta_1,\ldots,\theta_k\}$ is linearly
independent modulo $V_0=\{0\}$.  Now let $\pi_i: \Gg^* \into
\Gg^*/V_i$ be the natural projection, and suppose
$\{\pi_{i}(\theta_j)\ :\ \idx(\theta_j)>i\}$ is linearly
independent for all $i<m$.  Consider $X=\{\pi_m(\theta_j)\ :\
\idx(\theta_j)>m\}$. If there is a linear relation
$\pi_m(\theta_{p}) = \sum_{l\neq p} \ga_l\pi_m(\theta_l)$ among
these elements, then we may replace $\theta_{p}$ with
$\tilde{\theta}_p=\theta_p - \sum_{l\neq p}\ga_l \theta_l$, which
does not change $\Omega$ or affect linear independence modulo
$V_i,\ i<m$.   However, note that $\idx(\tilde\theta_p)\leq m$,
i.e. we have removed an element from $X$.  Reordering by degree
and repeating the argument, we may remove any linear relation
modulo $V_m$ in this way, proving the induction step.
\end{proof}

We require a simple linear algebra fact before moving on
to the theorem.
\begin{lem}\label{T:injection}
Let $V$ be a subspace of a vector space $W$, let
$\alpha\in\wedge^p V$, and suppose
$\{\theta_1,\ldots,\theta_m\}\subset W$ is linearly independent
modulo $V$. Then
$\alpha\wedge\theta_{1}\wedge\cdots\wedge\theta_{m}=0$ if and only
if $\alpha=0$.
\end{lem}
\begin{proof}
Let $\pi:W\rightarrow W/V$ be the projection and choose a
splitting $W\cong V\oplus W/V$;
$\alpha\wedge\theta_{1}\wedge\cdots\wedge\theta_{m}$ has a
component in $\wedge^p V\otimes\wedge^m(W/V)$ equal to
$\alpha\otimes \pi(\theta_1)\wedge\cdots\wedge\pi(\theta_m)$,
which vanishes if and only if $\alpha=0$.
\end{proof}

\begin{theorem}\label{pure spinor}  Any left-invariant generalized
complex structure on a nilmanifold must be generalized Calabi-Yau.
That is, any left-invariant global trivialization $\rho$ of the
canonical bundle must be a closed differential form.  In
particular, any left-invariant complex structure has
holomorphically trivial canonical bundle.
\end{theorem}
\begin{proof}
Let $\rho=e^{B+i\omega}\Omega$ be a left-invariant trivialization
of the canonical bundle such that
$\Omega=\theta_1\wedge\cdots\wedge\theta_k$, with
$\{\theta_1,\ldots,\theta_k\}$ ordered according to
Lemma~\ref{T:trivial}. Then the integrability condition
$d\rho=(X+\xi)\cdot\rho$ is equivalent to the condition
\begin{equation}\label{E:equi}
d(B+i\omega)\wedge\Omega + d\Omega =
(i_X(B+i\omega))\wedge\Omega+i_X\Omega+\xi\wedge\Omega.
\end{equation}
The degree $k+1$ part of~(\ref{E:equi}) states that
\begin{equation}\label{E:eq1}
d\Omega = i_X(B+i\om)\wedge\Omega + \xi\wedge\Omega.
\end{equation}
Taking wedge of \eqref{E:eq1} with $\theta_i$ we get
\begin{equation*}
d\theta_i \wedge \theta_1\wedge\cdots\wedge\theta_k=0,\ \ \forall
i.
\end{equation*}
Now, let $\{\theta_1,\ldots,\theta_j\}$ be the subset with
nilpotent degree $\leq \idx(d\theta_i)$.  Note that $j<i$ since
$\idx(\theta_i)=\idx(d\theta_i)+1$.  Then since
\begin{equation*}
(d\theta_i\wedge\theta_1\wedge\cdots\wedge\theta_j)\wedge\theta_{j+1}\wedge\cdots\wedge\theta_k
= 0,
\end{equation*}
we conclude from Lemma~\ref{T:injection} that
\begin{equation}\label{E:important}
d\theta_i\wedge\theta_1\wedge\cdots\wedge\theta_j=0, \mbox{ with } j<i.
\end{equation}
Since this argument holds for all $i$, we conclude that
$d\Omega=0$. The degree $k+3$ part of~(\ref{E:equi}) states that
$d(B+i\omega)\wedge\Omega=0$, and so we obtain that $d\rho =
e^{B+i\omega}d\Omega=0$, as required.
\end{proof}

Equation \eqref{E:important} shows that the integrability
condition satisfied by $\rho$ leads to constraints on
$\{\theta_1,\ldots,\theta_k\}$. Since these will be used
frequently in the search for \gcss, we single them out as follows.

\begin{cor}\label{T:ideal}
Assume $\{\theta_1,\ldots,\theta_k\}$ are chosen according to
Lemma \ref{T:trivial}. Then
\begin{equation}\label{idealtheta}
d\theta_i \in \mc{I}(\{\theta_j : \idx(\theta_j) <
\idx(\theta_i)\}),
\end{equation}
where $\mc{I}(~)$ denotes the
ideal generated by its arguments. Since $\nil(\theta_i)$ is weakly
increasing, we have, in particular,
\begin{equation*}
d\theta_i\in\mathcal{I}(\theta_1,\ldots,\theta_{i-1}).
\end{equation*}
\end{cor}

\begin{example}\label{T:simple case1}
Since $d\theta_1 \in \mc{I}(0)$, we see that $\theta_1$ is always
closed, and therefore it lies on $V_1$ or, equivalently, $\idx(\theta_1)=1$.
\end{example}

So far, we have described constraints deriving from the
integrability condition on $\rho$.  However, nondegeneracy (in
particular, $\Omega\wedge\overline{\Omega}\neq 0$) also places
constraints on the $\theta_i$ appearing in the decomposition of
$\gO$.  The following example illustrates this.

\begin{example}\label{T:simple case2}
If $\theta_1,\ldots,\theta_j\in V_i$, then nondegeneracy implies
that $\dim V_i \geq 2j$.  For a fixed nilpotent algebra, this
places an upper bound on the number of $\theta_j$ which can be
chosen from each $V_i$.
\end{example}

In the next lemma we prove a similar, but more subtle constraint
on the 1-forms~$\theta_i$.

\begin{lem}\label{T:highest index}
Assume that $\{\theta_1,\ldots,\theta_k\}$ are chosen according to
Lemma \ref{T:trivial}. Suppose that no $\theta_i$ satisfies
$\idx(\theta_i)=j$, but that there exists $\theta_l$ with
$\idx(\theta_l) =j+1$. Then $\theta_l\wedge\overline\theta_l\neq
0$ modulo $V_j$ (i.e. in $\wedge^2(V_{j+1}/V_j)$), and in
particular $V_{j+1}/V_j$ must have dimension 2 or greater.
\end{lem}

\begin{proof}
Assume that the hypotheses hold but
$\theta_l\wedge\overline\theta_l=0$ modulo $V_j$.  Because of
this, it is possible to decompose $\theta_l=v+\alpha$, where
$\idx(\alpha)< j+1$, and such that $v\wedge\overline v = 0$.
Therefore, up to multiplication of $\theta_l$ (and therefore
$\Omega$) by a constant, $v$ is real.

By Corollary \ref{T:ideal}, $d\theta_l\in
\mathcal{I}(\{\theta_i : \nil(\theta_i)< j+1\})$. By hypothesis, there is no $\theta_i$ with nilpotent degree $j$, therefore we obtain
\begin{equation*}
dv +d\ga = \sum_{\idx(\theta_i) < j} \xi_i \wedge \theta_i,
\end{equation*}
where $\xi_i\in\Gg^*$.   Since $\idx(dv)=j$, there is an element
$w \in \Gg^{j-1}$ such that $i_w dv \neq 0$.  On the other hand,
the nilpotent degrees of $d\ga$ and the $\{\theta_i\}$ in the sum above are
less than $j$, hence interior product with $w$ annihilates all
these forms.  In particular,
\begin{equation*} 0 \neq i_w dv = \sum_{\idx(\theta_i) < j}(i_w\xi_i)\theta_i.
\end{equation*}
Therefore, $i_w\xi_i$ is nonzero for some $i$. But, the left hand side
is real, thus
\begin{equation*}0 = i_w dv \wedge \overline{i_w dv}
= \left(\sum_{\idx(\theta_i) < j}(i_w\xi_i)\theta_i\right) \wedge
\left(\sum_{\idx(\theta_i) <
j}(i_w\overline{\xi_i})\overline{\theta_i} \right).
\end{equation*}
By the nondegeneracy condition, the right hand side is nonzero,
which is a contradiction.  Hence
$\theta_l\wedge\overline\theta_l\neq 0$ modulo $V_j$.
\end{proof}

From this lemma, we see that if $\dim V_{j+1}/V_j = 1$ occurs in a
nilpotent Lie algebra, then it must be the case, either that some
$\theta_i$ has nilpotent degree $j$, or that no $\theta_i$ has
nilpotent degree $j+1$.  In this way, we see that the size of the
nilpotent steps $\dim V_{j+1}/V_j$ may constrain the possible
types of left-invariant generalized complex structures, as we now
make precise.

\begin{theorem}\label{T:maximal}
Let $M^{2n}$ be a nilmanifold with associated Lie algebra
$\frak{g}$. Suppose there exists a $j>0$ such that, for all $i\geq
j$,
\begin{equation}\label{jump}
\dim \left({V_{i+1}}/{V_{i}}\right) = 1.
\end{equation}
Then $M$ cannot admit left-invariant generalized complex
structures of type $k$ for $k \geq 2n - \nil(\Gg) +j $.

In particular, if $M$ has maximal nilpotency index (i.e. $\dim
V_1=2,\ \dim V_i/V_{i-1}=1\ \forall i>1$), then it does not admit
\gcss\ of type $k$ for $k\geq 2$.
\end{theorem}
\begin{proof}

First observe that for any nilpotent Lie algebra $\nil(\Gg)\leq
2n-1$, so the theorem only restricts the existence of structures
of type $k>1$.

According to Lemma~\ref{T:highest index}, if none of the
$\theta_i$ have nilpotent degree $j$, then there can be none with
nilpotent degree $j+1, j+2, \ldots$ by the condition~(\ref{jump}).
Hence we obtain upper bounds for the nilpotent degrees of
$\{\theta_1,\ldots,\theta_k\}$, as follows.  First, $\theta_1$ has
nilpotent degree 1 (by Example~\ref{T:simple case1}).  Then, if
$\idx(\theta_2)\geq j+2$, this would imply that no $\theta_i$ had
nilpotent degree $j+1$, which is a contradiction by the previous
paragraph.  Hence $\idx(\theta_2) < j+2$.  In general,
$\idx(\theta_i) < j+i$.  Suppose that $M$ admits a generalized
complex structure of type $k>1$.  Then we see that
$\idx(\theta_k)< j+k$. By Example~\ref{T:simple case2}, we see
this means that $\dim V_{j+k-1}\geq 2k$.

On the other hand, $\dim V_{j+k-1} = 2n - \dim \Gg^*/V_{j+k-1}$,
and since $\Gg^*=V_{\nil(\Gg)}$, we have
\begin{equation*}
\Gg^*/V_{j+k-1} = V_{\nil(\Gg)}/V_{j+k-1} \cong
V_{\nil(\Gg)}/V_{\nil(\Gg)-1}\oplus\cdots\oplus V_{j+k}/V_{j+k-1},
\end{equation*}
and the $\nil(\Gg)-j-k+1$ summands on the right have dimension 1,
by hypothesis.  Hence $\dim V_{j+k-1} = 2n-\nil(\Gg)+j+k-1$, and
combining with the previous inequality, we obtain
\begin{equation*}
k < 2n-\nil(\Gg)+j,
\end{equation*}
as required. For the last claim, observe that $M^{2n}$ has maximal
nilpotency index when $\nil(M) = 2n-1$, in which case
$\eqref{jump}$ holds for $j=1$.
\end{proof}

Constraints beyond those mentioned here may be obtained if one
considers the fact that $\gO\wedge \overline{\gO}$ defines a
foliation and that $\go$ restricts to a symplectic form on each
leaf.  Both the leafwise nondegeneracy of $\omega$ and the
requirement of being closed on the leaves lead to useful
constraints on what types of generalized complex structure may
exist, as we shall see in the following sections.

\section{Generalized complex structures on
6-nilmanifolds}\label{six}

In this section, we turn to the particular case of 6-dimensional
nilmanifolds. The problem of classifying those which admit
left-invariant complex (type 3) and symplectic (type 0) structures
has already been solved \cite{Sa,GK}, so we are left with the task
of determining which 6-nilmanifolds admit left-invariant
generalized complex structures of types $1$ and $2$. The result of
our classification is presented in Table~1, where explicit
examples of all types of left-invariant generalized complex
structures are given, whenever they exist.  The main results
establishing this classification are Theorems \ref{T:2,1} and
\ref{T:1,2}.  Throughout this section we often require the use of
a Malcev basis $\{e_1,\ldots,e_6\}$ as well as its dual basis
$\{\del_1,\ldots,\del_6\}$.  We use $e_{i_1\cdots i_p}$ as an
abbreviation for $e_{i_1}\wedge\cdots\wedge e_{i_p}$.

\subsection{Generalized complex structures of type 2}

By the results of the last section, a left-invariant structure of
type $2$ is given by a closed form $\rho = \exp(B+i\om) \theta_1
\theta_2$, where
$\om\wedge\theta_1\theta_2\overline{\theta}_1\overline{\theta}_2\neq
0$. As a consequence of Theorem \ref{T:maximal}, any 6-nilmanifold
with maximal nilpotence step cannot admit this kind of structure.
We now rule out some additional nilmanifolds, using Lemma
\ref{T:highest index}.

\begin{lem}\label{lemmafirst}
If a 6-nilmanifold $M$ has nilpotent Lie algebra given by
$(0,0,0,12,14,-)$, and has nilpotency index $4$, then $M$ does not
admit left-invariant generalized complex structures of type $2$.
\end{lem}
\begin{proof}
Suppose that $M$ admits a structure $\rho$ of type $2$.  Since $M$
has nilpotency index $4$, $V_{i+1}/V_i$ is 1-dimensional for $i
\geq 1$. From $d\theta_1=0$ and Lemma \ref{T:highest index}, we
see that $\theta_1= z_1e_1+z_2e_2+z_3e_3$ and $\nil(\theta_2) \leq
2$, thus $\theta_2= w_1e_1+w_2e_2+w_3e_3+w_4e_4$. The conditions
$d(\theta_1 \theta_2)=0$ and
$\theta_1\theta_2\overline{\theta_1}\overline{\theta_2} \neq 0$
together imply $z_3=0$.  Further, the annihilator of
$\theta_1\theta_2\overline{\theta_1}\overline{\theta_2}$ is
generated by $\{\del_5,\del_6\}$.  Hence, the nondegeneracy
condition $\omega^2\wedge\Omega\wedge\overline{\Omega}\neq 0$
implies that
$$B+i\om = (k_1e_1 +
\dots k_4 e_4 + k_5e_5)e_6 + \alpha,$$ where $k_5 \neq 0$ and $\ga
\in \wedge^2\left<e_1,\cdots,e_5\right>$. But, using the structure
constants, we see that $d\rho$ must contain a nonzero multiple of
$e_6$, and so cannot be closed.
\end{proof}

\begin{lem}\label{lemmasecond}
Nilmanifolds associated to the algebras defined by
\begin{align*}
&(0,0,0,12,14,13-24),\\
&(0,0,0,12,14,23+24)
\end{align*}
do not admit left-invariant generalized complex structures of type
$2$.
\end{lem}
\begin{proof}
Each of these nilmanifolds has $\nil(\Gg)=3$, with $\dim V_1=3$
and $\dim V_2/V_1=1$.  Suppose either nilmanifold admitted a
structure $\rho$ of type $2$.  If $\nil(\theta_2) = 2$, we could
use the argument of the previous lemma to obtain a contradiction.
Hence, suppose $\idx(\theta_2)=3$.  Lemma \ref{T:highest index}
then implies that $\theta_2\wedge\overline{\theta_2}\neq 0\pmod
{V_2}$, which means it must have nonzero $e_5$ and $e_6$
components.

Now, if $\theta_1$ had a nonzero $e_3$ component, then
$d\theta_2\wedge\theta_1$ would have nonzero $e_{234}$ component,
contradicting~(\ref{idealtheta}).  Hence
\begin{equation}
\theta_1= z_1 e_1 + z_2 e_2 \qquad \theta_2= \sum_{i=1}^6 w_i e_i.
\end{equation}
But for these, the coefficient of $e_{123}$ in
$d\theta_2\wedge\theta_1$ would be $- z_2 w_6$ (for the first
nilmanifold) or $z_1 w_6$ (for the second), in each case implying
$\theta_1\wedge\overline{\theta_1}=0$, a contradiction.
\end{proof}

\begin{lem}\label{T:lemma3}
Nilmanifolds associated to the algebras defined by
\begin{align*}
&(0,0,12,13,23,14),\\
&(0,0,12,13,23,14-25)
\end{align*}
do not admit left-invariant generalized complex structures of type
$2$.
\end{lem}
\begin{proof}
Each of these nilmanifolds have $\nil(\Gg)=4$, with $\dim V_1=2$,
$\dim V_2=3$, and $\dim V_3=5$.  Suppose either nilmanifold
admitted a structure $\rho$ of type $2$. $V_4/V_3$ is
1-dimensional and so Lemma \ref{T:highest index} implies that
$\idx(\theta_2) \neq 4$. Since $\theta_1$ is closed and satisfies
$\theta_1\overline{\theta}_1\neq 0$, we may rescale it to obtain
$\theta_1=e_1+z_2e_2$. The condition $d\theta_2 \in
\mathcal{I}(\theta_1)$ implies we can write $\theta_2 = w_2e_2 +
w_3e_3 + w_4(e_4 + z_2e_5)$. Now let
\begin{equation*}
B+i\om=\sum_{i<j}k_{ij}e_{ij}
\end{equation*}
and differentiate $\rho$ using the structure constants.  In both
cases, $d\rho=0$ implies $B+i\om = \xi_1\theta_1 + \xi_2\theta_2$
for 1-forms $\xi_i$. Therefore \om\ is degenerate on the leaves
defined by
$\Ann(\theta_1\theta_2\overline{\theta_1}\overline{\theta_2})$,
contradicting the requirement
$\omega^2\wedge\Omega\wedge\overline\Omega\neq 0$.
\end{proof}

\begin{theo}{2,1}
The only 6-dimensional nilmanifolds not admitting left-invariant
generalized complex structures of type 2 are those with maximal nilpotency index and those excluded by
Lemmas~\ref{lemmafirst}, \ref{lemmasecond}, and \ref{T:lemma3}.
\end{theo}
\begin{proof}
In Table 1, we provide explicit left-invariant generalized complex
structures of type 2 for all those not excluded by the preceding
lemmas and Theorem \ref{T:maximal}.
\end{proof}

\subsection{Generalized complex structures of type 1}

A left-invariant generalized complex structure of type $1$ is
given by $\rho =\exp(B+i\om)\theta_1$, where $\om^2\wedge
\theta_1\overline{\theta_1} \neq 0$.  Note that this implies that
$\omega$ is a symplectic form on the $4$-dimensional leaves of the
foliation determined by $\theta_1\wedge\overline{\theta_1}$.

\begin{theo}{1,2}
The only 6-nilmanifolds which do not admit left-invariant
generalized complex structures of type $1$ are those associated to
the algebras defined by $(0,0,12,13,23,14)$ and
$(0,0,12,13,23,14-25)$.
\end{theo}
\begin{proof}
In Table 1, we provide explicit forms defining type 1 structures
for all 6-nilmanifolds except the two mentioned above.

Suppose the nilmanifold is associated to the Lie algebra
$(0,0,12,13,23,14)$. Then up to a choice of Malcev basis, a
generalized complex structure of type 1 can be written
\begin{equation*}
\rho = \exp(B+i\om)(e_1+z_2e_2), \qquad B+i\om =
\sum_{i<j}k_{ij}e_{ij}.
\end{equation*}
The condition $d\rho=0$ implies that
\begin{equation*}
(-k_{36}e_{314} +k_{45}e_{135} - k_{45}e_{423} + k_{46}e_{136} +
k_{56}e_{236}-k_{56}e_{514})(e_1+z_2e_2)=0.
\end{equation*}
The vanishing of the $e_{1245}$, $e_{1236}$, $e_{1235}$, and
$e_{1234}$ components imply successively that $k_{56}$, $k_{46}$,
$k_{45}$, and $k_{36}$ all vanish.

The leaves of the distribution $\Ann(\theta_1\overline{\theta_1})$
are the tori generated by $\del_3,\del_4,\del_5,\del_6$, and the
previous conditions on $B+i\omega$ imply that on these leaves,
\om\ restricts to $\mbox{Im}(k_{34})e_{34} +
\mbox{Im}(k_{35})e_{35}$ which is degenerate, contradicting
$\om^2\wedge\theta_1\overline{\theta}_1\neq 0$. An identical
argument holds for the nilpotent algebra $(0,0,12,13,23,14-25)$.
\end{proof}

\section{$\beta$-transforms of generalized complex
structures}\label{iwa}

In this section, we will use Proposition~\ref{T:beta field} to
show that any left-invariant complex structure on a nilmanifold
may be deformed into a left-invariant generalized complex
structure of type 1.  By connecting the type 3 and type 1
structures in this way, we go on to show that the two disconnected
components of the left-invariant complex moduli space on the
Iwasawa manifold may be joined by paths of generalized complex
structures.

\begin{theo}{complex}
Every left-invariant complex structure  $2n$-dimensional
nilmanifold can be deformed, via a \gb-field, into a
left-invariant generalized complex structure of type $n-2$.
\end{theo}
\begin{proof}
According to Proposition \ref{T:beta field}, such a deformation
can be obtained if we find a holomorphic Poisson structure.  Let
us construct such a bivector $\beta$.

By Theorem~\ref{pure spinor}, a left-invariant complex structure
on a nilmanifold has a holomorphically trivial canonical bundle.
Let $\Omega=\theta_1\wedge\cdots\wedge\theta_n$ be a holomorphic
volume form, and choose the $\theta_i$ according to
Lemma~\ref{T:trivial}, so that, by Corollary~\ref{T:ideal}, the
differential forms
$\theta_1,\theta_1\wedge\theta_2,\ldots,\theta_1\wedge\cdots\wedge\theta_n$
are all holomorphic.  Now let $\{x_1,\ldots,x_n\}$ be a dual basis
for the holomorphic tangent bundle.  By interior product with
$\Omega$, we see that the multivectors $x_n,x_{n-1}\wedge
x_n,\ldots,x_1\wedge\cdots\wedge x_n$ are all holomorphic as well.
In particular, we have a holomorphic bivector $\beta=x_{n-1}\wedge
x_n$.  Calculating the Schouten bracket of this bivector with
itself, we obtain
\begin{equation*}
[\beta,\beta]=[x_{n-1}\wedge x_n,x_{n-1}\wedge
x_n]=2[x_{n-1},x_n]\wedge x_{n-1}\wedge x_n=0,
\end{equation*}
where the last equality follows from the fact that
$[x_{n-1},x_n]\in\left<x_{n-1},x_n\right>$, since
$\theta_i([x_{n-1},x_n])=-d\theta_i(x_{n-1},x_n)=0$ for $i<n-1$,
by Corollary~\ref{T:ideal}.

Hence $\beta$ gives rise to a deformation of the generalized
complex structure. The resulting structure $\tilde\rho$ is given
by the following differential form:
\begin{equation*}
\tilde\rho=e^{i_\beta}\rho = \rho + i_\beta\rho =
e^{\theta_{n-1}\wedge\theta_n}\theta_1 \wedge \dots \wedge
\theta_{n-2},
\end{equation*}
and we see immediately that it is a left-invariant generalized
complex structure of type $n-2$.
\end{proof}

In \cite{KS}, Ketsetzis and Salamon study the space of
left-invariant complex structures on the Iwasawa nilmanifold. This
manifold is the quotient of the complex 3-dimensional Heisenberg
group of unipotent matrices by the lattice of unipotent matrices
with Gaussian integer entries. As a nilmanifold, it has structure
$(0,0,0,0,13-24,14+23)$. Ketsetzis and Salamon observe that the
space of left-invariant complex structures with fixed orientation
has two connected components which are distinguished by the
orientation they induce on the complex subspace $\langle \del_1,
\del_2,\del_3, \del_4 \rangle$.

\vskip6pt \noindent {\bf Connecting the two components.} Consider
the left-invariant complex structures defined by the closed
differential forms $\rho_1 = (e_1+ie_2)(e_3+ie_4)(e_5+ie_6)$ and
$\rho_2 = (e_1+ie_2)(e_3-ie_4)(e_5-ie_6)$.  These complex
structures clearly induce opposite orientations on the space
$\langle \del_1, \dots, \del_4 \rangle$, so lie in different
connected components of the moduli space of left-invariant complex
structures.

By Theorem \ref{T:complex}, the first complex structure can be
deformed, by the $\beta$-field $\beta_1 = \frac{-1}{4}(x_3 -
ix_4)(x_5-ix_6)$ into the generalized complex structure
$$e^\beta \rho_1 = e^{-(e_{35} - e_{46}) - i(e_{45} + e_{36})}(e_1+ie_2),$$
and then, by the action of the closed $B$-field $B_1= e_{35} - e_{46}$, into
$$\rho = e^{ - i(e_{45} + e_{36})}(e_1+ie_2).$$

On the other hand, the second complex structure can be deformed,
via the $\beta$-field $\beta_2 = \frac{1}{4}(x_3 +
ix_4)(x_5+ix_6)$, into the type 1 generalized complex structure
$$e^\beta \rho_2 = e^{(e_{35} - e_{46}) - i(e_{45} + e_{36})}(e_1+ie_2),$$
and then by the $B$-field $B_2 = -(e_{35} - e_{46})$ into
$$ \rho =e^{ - i(e_{45} + e_{36})}(e_1+ie_2),$$
which is the same generalized complex structure obtained from
$\rho_1$.

Therefore, by deforming by $\beta$- and $B$-fields, the two
disconnected components of the moduli space of left-invariant
complex structures can be connected, through generalized complex
structures.

\section{An 8-dimensional nilmanifold}\label{8d}

We have established that all 6-dimensional nilmanifolds admit
generalized complex structures.  One might ask whether every
even-dimensional nilmanifold admits left-invariant generalized
complex geometry.  In this section, we answer this question in the
negative, by presenting an 8-dimensional nilmanifold which does
not admit any type of left-invariant generalized complex
structure.

\begin{ex}{8d}
Consider a nilmanifold $M$ associated to the 8-dimensional
nilpotent Lie algebra defined by
\begin{equation*}
(0,0,12,13,14,15,16,36-45-27).
\end{equation*}

Since it has maximal nilpotency index, Theorem \ref{T:maximal}
implies that it may only admit left-invariant generalized complex
structures of types 1 and 0.  We exclude each case in turn:

\begin{itemize}
\item \emph{Type 1}: Suppose there is a type 1 structure, defined
by the left-invariant form $\rho = e^{B+i\om}\theta_1$.  Then
$d\theta_1$ =0 and $\theta_1\overline{\theta}_1\neq 0$ imply that
$\theta_1\overline{\theta}_1$ is a multiple of $e_{12}$ and
therefore $\om$ must be symplectic along the leaves defined by
$\langle\del_3,\ldots,\del_8\rangle$.  These leaves are actually
nilmanifolds associated to the nilpotent algebra defined by
$(0,0,0,0,0,12+34)$, which admits no symplectic structure, and so
we obtain a contradiction.
\item \emph{Type 0}:
The real second cohomology of $M$ is given by
\begin{equation*}
H^2(M,\RR) = \langle e_{23},e_{34}-e_{25},e_{17}\rangle,
\end{equation*}
and since $e_8$ does not appear in any of its generators, it is
clear that any element in $H^2(M,\RR)$ has vanishing fourth power,
hence excluding the existence of a symplectic structure.
\end{itemize}
In this way, we see that the 8-dimensional nilmanifold $M$ given
above admits no left-invariant generalized complex structures at
all.
\end{ex}

\setlength{\oddsidemargin}{0 cm} \setlength{\evensidemargin}{-0.25
cm} \setlength{\topmargin}{2.0 cm} \setlength{\footskip}{-4.0cm}
\setlength{\parindent}{2em} \setlength{\textwidth}{17.5 cm}
\setlength{\textheight}{22 cm} \setlength{\baselineskip}{40pt}
\pagestyle{empty}
\begin{landscape}
\begin{center}
{\tiny
\begin{tabular}{|l c c l l l l|}
\hline
Nilmanifold class & $b_1$ & $b_2$ & Complex (type 3) &Type 2 & Type 1 & Symplectic (type 0)\\
\hline $(0,0,12,13,14,15)$ & 2 & 3 &  -- & -- &
$(1+i2)\exp{i(36-45)}$ & $16 + 34
- 25$ \\
$(0,0,12,13,14,34+52)$ & 2 & 2 & -- & -- & $(1+i2)\exp(- 45 +  36 + i(36+ 45))$ & --\\
$(0,0,12,13,14,23+15)$ & 2 & 3 & -- & -- & $(1+i2)\exp i(36-45)$ &
$16 +
24 + 34 -25$\\
$(0,0,12,13,23,14)$ & 2 & 4 & -- &  -- & -- & $15 + 24 + 34 - 26$\\
$(0,0,12,13,23,14-25)$ & 2 & 4 & -- & -- & -- & $15 + 24 - 35 + 16$\\
$(0,0,12,13,23,14+25)$ & 2 & 4 & $(1+i2)(4+i5)(3+i6)$ & $(1 +
i2)(4+i5)\exp i(36)$ &
$(1 + i2)\exp(43-56+i(46 - 35))$ & $15 + 24 + 35 + 16$\\
$(0,0,12,13,14+23,34+52)$ & 2 & 2 & -- & -- & $(1+i2)\exp(45-35+36+i(-36+45-16))$  & --\\
$(0,0,12,13,14+23,24+15)$ & 2 & 3 & -- & -- & $(1+i2)\exp(2 \times 35 +
i(36- 45))$ & $16 + 2 \times 34 - 25$\\
\hline $(0,0,0,12,13,14+35)$ & 3 & 5 & -- & $(1+2+i3)(4+i5)\exp
i(26) $ &$(1+i2)\exp i(
 36 + 45) $ & --\\

$(0,0,0,12,13,14+23)$ & 3 & 6 & $(1+i2)(3 - 2 \times i4)(5+
2\times i6)$ & $(2+i3)(4+i5)\exp i(16-34) $ & $(1+i2)\exp
i(36+45)$ & $16 -2\times
34 - 25$ \\

$(0,0,0,12,13,24)$ & 3 & 6 & $(1+i2)(3+4+i4)(5+6-i6)$ &
$(1+2+i3)(4+i5)\exp i(26) $ &
 $(1+i2)\exp i(35+46)$ & $26+14+35 $ \\

$(0,0,0,12,13,14)$ & 3 & 6 & $(1+i2)(3+i4)(5+i6)$ & $(2+i3) (4+i5)
 \exp i(16)$ & $(1+i2)\exp(35-46+i(36+45))$ & $16+24+35$  \\

$(0,0,0,12,13,23)$ & 3 & 8 & $(1+i2)(3+i4)(5+i6)$ &
$(1+i2)(5+i6)\exp i(16-34)$ & $(1+i2)\exp(35-46 + i(36 + 45))$ &
$15 + 24
+36$ \\

$(0,0,0,12,14,15+23)$ & 3 & 5 & -- & -- & $(1+i3)\exp i(26-45)$ & $13 +26
- 45$ \\

$(0,0,0,12,14,15 + 23 + 24)$ & 3 & 5 & -- & -- & $(1+i3) \exp i(26-45)$ &
$13+26-45$\\

$(0,0,0,12,14,15+24)$ & 3 & 5 & -- & -- & $(1+i3) \exp i(26-45)$ &
$13+26-45$\\

$(0,0,0,12,14,15)$ & 3 & 5 & -- & -- & $(1+i3) \exp i(26-45)$ &
$13+26-45$\\

$(0,0,0,12,14,24)$ & 3 & 5 & $(1+i2)(3+i4)(5+i6)$ &
$(1+i2)(5+i6)\exp i(34)$ & $(1+i2)\exp(35-46 + i(45+36))$ & --\\

$(0,0,0,12,14,13+42)$ & 3 & 5 & $(1+i2)(3+i4)(2 \times 5 - i6)$ &
-- &  $(1+i2)\exp i(35 + 46)$ & $15 +
26 + 34$\\

$(0,0,0,12,14,23+24)$ & 3 & 5 & $(1+i2)(3+4+i3)(5+6+i6)$ & -- &
$(1+i2)\exp(35+46+i(35-46))$ & $16-34+25$\\

$(0,0,0,12,23,14+35)$ & 3 & 5 & -- & $(1+2+i3)(5+i4)\exp(3+i1)6$ & $(1+i2)\exp(36 + 45 +i(36 - 45))$ & --\\

$(0,0,0,12,23,14-35)$ & 3 & 5 & $(1+i3)(4-i5)(2+i6)$ &
$(1+i3)(4-i5)\exp i(26)$ & $(1+i3)\exp(24 + 56 + i(25+46))$ & --\\

$(0,0,0,12,14-23,15+34)$ & 3 & 4 & -- &  $(1+i2)(3+i4)\exp i(56)$
&
 $(2+i3) \exp i(16 + 35 + 45-26)$  & $16+35+24$ \\

$(0,0,0,12,14+23,13+42)$ & 3 & 5 & $(1+i2)(3-i4)(5+i6)$ &
$(1+i2)(3-i4)\exp i(56)$ & $(1+i2)\exp(35+46 + i(36-45))$ &
 $15+2\times 26+34$\\

\hline

$(0,0,0,0,12,15+34)$ & 4 & 6 & -- & $(1+i2)(3+i4)\exp i(56)$ & $(3+i4)\exp i(25+16)$ & --\\
$(0,0,0,0,12,15)$ & 4 & 7 & -- & $(1+i2)(3+i4)\exp i(56)$ &
$(1+i2)\exp i(34 + 56)$ & $16+25+34$\\
$(0,0,0,0,12,14+25)$ & 4 & 7 & $(1+i2)(4+i5)(3+i6)$ &
$(1+i2)(4+i5)\exp i(36)$ & $ (1+i2) \exp(34 + 56 + i(35 - 46))$& $13+26+45$\\
$(0,0,0,0,12,14+23)$ & 4 & 8 & $ (1+i2)(3-i4)(5+i6)$ &
$(1+i2)(3-i4)\exp i(56)$ & $(1+i2)\exp(35 + 46 +i(36 - 45))$ & $13+26+45$\\
$(0,0,0,0,12,34)$ & 4 & 8 & $(1+i2)(3+i4)(5+i6)$&
$(1+i2)(3+i4)\exp i(56)$ & $(1+i2)\exp(35-46 + i(45 + 36))$ & $15+36+24$\\
$(0,0,0,0,12,13)$ & 4 & 9 & $(2+i3)(1+i4)(5+i6)$ &
$(2+i3)(5+i6)\exp i(14)$ & $(2+i3)\exp(15-46 +i(16+45))$ & $16+25+34$ \\
($0,0,0,0,13+42,14+23)$ & 4 & 8 & $(1+i2)(3-i4)(5+i6)$ &
$(1+i2)(3+i4)\exp i(56)$ & $(1+i2)\exp(35+46+i(36-45))$ & $16+25+34$ \\
\hline $(0,0,0,0,0,12+34)$ & 5 & 9 & $(1+i2)(3+i4)(5+i6)$ &
$(1+i2)(3+i4)\exp i(56)$ &
$(1+i2)\exp i(36+45)$ & --\\
$(0,0,0,0,0,12)$ & 5 & 11 & $(1+i2)(3+i4)(5+i6)$ &
$(1+i2)(3+i4)\exp i(56)$ & $(1+i2)\exp i(36+45)$& $16+23+45$\\
 \hline
$(0,0,0,0,0,0)$ & 6 & 15 & $(1+i2)(3+i4)(5+i6)$ &
$(1+i2)(3+i4)\exp i(56)$ &
$(1+i2)\exp i(36+45)$ & $12+34+56$\\
\hline
\end{tabular}\label{table}}
\vskip18pt {\small {\it Table 1: Differential forms defining
left-invariant Generalized Calabi-Yau structures. The symbol `-'
denotes nonexistence.}}
\end{center}
\end{landscape}

\end{document}